\documentclass{amsart}
\usepackage{amscd,amssymb,subfigure}
\usepackage[all,dvips]{xy}
\theoremstyle{plain}
\newtheorem{lemma}{Lemma}
\newtheorem{theorem}[lemma]{Theorem}

\newtheorem*{thmplus}{Theorem $\mathbf 1^+$}
\newtheorem{corollary}[lemma]{Corollary}
\theoremstyle{definition}
\newtheorem{ex}{Example}
\newtheorem{remark}{Remark}
\def\ba{{\Delta}} % base
\def\affphi{{\widetilde{\Phi}}} % affine root system
\def\ca{{\mathcal A}}
\def\car{{\ca^{H_0}}}
\def\cax{{\ca^{H_x}}}
\def\set#1{\left\{#1\right\}}
\def\abs#1{\left\vert#1\right\vert}
\def\R{{\mathbb R}}
\newxyColor{unbd}{0.74 0.35 0.92}{rgb}{}
\subjclass[2000]{Primary
20F55,  %% Reflection and Coxeter groups
Secondary
52C35%  %% Arrangements of points, flats, hyperplanes
}
\keywords{Root systems, real hyperplane arrangements, reflection groups}

\title{A note on de Concini and Procesi's curious identity}

\author[G. Denham]{Graham Denham}
\address{Department of Mathematics, University of Western Ontario,
London, ON  N6A 5B7}
\urladdr{{http://www.math.uwo.ca/\~{}gdenham}}
\thanks{Partially supported by a grant from NSERC of Canada.}
\begin{document}
\begin{abstract}
We give a short, case-free and combinatorial proof of de Concini
and Procesi's formula from \cite{dCP06} 
for the volume of the simplicial cone spanned
by the simple roots of any finite root system.  The argument presented
here also extends their formula to include 
the non-crystallographic root systems.
\end{abstract}
\maketitle
\section{Introduction}
Let $\Phi\subseteq\R^n$ be a finite root system with
base $\ba$, and let $W=W(\Phi)$ denote the reflection group of $\Phi$.
Let $\sigma_\ba$ be the positive cone spanned by the set of 
simple roots $\Delta$:
\begin{equation}\label{def:sigma}
\sigma_\ba=\set{\sum_{\alpha\in\ba}c_\alpha\alpha\colon
c_\alpha\in\R_{>0}\hbox{~for all $\alpha\in \ba$}}.
\end{equation}
Let $C_\ba$ be the normal cone to $\sigma_\ba$: this is usually called
the fundamental chamber in the arrangement $\ca$ of reflecting hyperplanes
of $W$.  
If $\tau$ is a cone in $\R^n$, define the volume of $\tau$ as
$\nu(\tau)={\rm vol}(\tau\cap D^n)/
{\rm vol}\;D^n$, where $D^n$ is the unit ball centered at the origin.
Finally, let $\set{d_1, d_2,\ldots, d_n}$ denote the degrees of
$W$: we refer to \cite{humpcoxbook} for background and notation.

Recall that the action of $W$ on $\R^n$ by reflections is free 
on the complement of the hyperplanes $\ca$.  The induced action on 
chambers is simply transitive.
Since the chambers partition the complement of $\ca$
and $W$ acts by isometries,
$\nu(gC_{\ba})=1/\abs{W}=1/\prod_{i=1}^n d_i,
$
for any chamber $g C_\ba$.

While not so straightforward, it turns out that the volume of the 
cone $\sigma_\Delta$ is also rational, and has a nice expression:
\begin{theorem}[Theorem~1.3 in~\cite{dCP06}]\label{th:vol}
If $\Phi$ is crystallographic, the volume of the cone $\sigma_\Delta$
is
\begin{equation}\label{eq:volformula}
\nu(\sigma_\ba)=\prod_{i=1}^n \frac{d_i-1}{d_i}.
\end{equation}
\end{theorem}
De Concini and Procesi derive this result from the ``curious identity''
of their title.  Their proof of the identity is accompanied by a note
by Stembridge that gives an elegant, alternate proof via character
theory.

The purpose of this note is to offer yet another argument.  Using the
combinatorial theory of real hyperplane arrangements, one 
can prove \eqref{eq:volformula} directly, in slightly more generality
(\S\ref{sec:vol}).  
Then, in the crystallographic case,
de Concini and Procesi's identity is recovered by adding up normal cones
around the fundamental alcove of the associated affine root system
$\affphi$ (in \S\ref{sec:id}).

\section{The volume formula}\label{sec:vol}
Let $V\subseteq\R^n$ consist of the union of the reflecting hyperplanes,
together with those vectors in the span of any proper subset of any base
$g\ba$.  Clearly $\R^n-V$ is a dense, open subset of $\R^n$.  The key
result is the following, whose proof appears at the end of this
section.

\begin{theorem}\label{th:count}
For any $x\in\R^n-V$, the number of $g\in W$ for which $x\in g\sigma_{\ba}$
is independent of $x$ and equal to $\prod_{i=1}^n(d_i-1)$.
\end{theorem}
In another formulation,
\begin{corollary}
For a finite root system $\Phi$ and $x\in\R^n-V$, the number of choices
of base $\Delta$ for $\Phi$ for which $x$ is in the positive cone of $\Delta$
equals $\prod_{i=1}^n(d_i-1)$.
\end{corollary}
\begin{proof}
If $\ba$, $\ba'$ are both bases for $\Phi$, then $\ba'=g\ba$ for some
$g\in W$, and $\sigma_{\ba'}=g\sigma_\ba$.
\end{proof}
Since each cone $g\sigma_\ba$ has the same volume, 
\begin{eqnarray*}
\abs{W}\cdot \nu(\sigma_\ba)&=&\sum_{g\in W}\nu(g\sigma_\ba)\\
&=&\prod_{i=1}^n (d_i-1)
\end{eqnarray*}
by Theorem~\ref{th:count}, and we obtain the volume formula
as a corollary:
\begin{thmplus}\label{th:main}
If $\Phi$ is any finite root system, the volume of the cone $\sigma_\ba$ is
$$
\nu(\sigma_\ba)=\prod_{i=1}^n \frac{d_i-1}{d_i}.
$$
\end{thmplus}
(Note that, if the rank of $\Phi$ is less than $n$, the least degree
is $1$, and both sides are zero.)
\subsection{Hyperplane arrangements}
The terminology used below may be found in the book of Orlik
and Terao~\cite{ot}.  We recall a collection of
hyperplanes $\ca$ in $\R^n$ is {\em central} if all $H\in\ca$ contain the
origin, and {\em essential} if the collection of normal vectors span
$\R^n$.  

Recall that $\ca$ has an intersection lattice $L(\ca)$ of
subspaces, ranked
by codimension.  The Poincar\'e polynomial of $\ca$ is defined to be
$$
\pi(\ca,t)=\sum_{X\in L(\ca)}\mu(\widehat{0},X)(-t)^{{\rm rank}(X)},
$$
where $\mu$ is the M\"obius function.  If $\ca$ is essential, $\pi(\ca,t)$
is a polynomial of degree $n$.  The following classical theorem 
is a main ingredient in our proof.
\begin{theorem}[\cite{os80b}]
If $\ca=\ca(\Phi)$ is an arrangement of (real) reflecting hyperplanes, 
then
\begin{equation}\label{eq:terao}
\pi(\ca,t)=\prod_{i=1}^n(1+(d_i-1)t),
\end{equation}
where $\set{d_i}$ are the degrees of the reflection group.
\end{theorem}
If $H_0$ is any hyperplane (not necessarily through the origin), let $\car$
denote the set $\set{H\cap H_0:H\in\ca}$, regarded as a hyperplane arrangement
in $H_0$.  We say $H_0$ is in general position to $\ca$ if 
$X\cap H_0$ is nonempty for all nonzero subspaces $X\in L(\ca)$.
\begin{lemma}\label{lem:trunc}
If $H_0$ is in general position to a central arrangement $\ca$ in $\R^n$, 
then the 
number of bounded chambers in $\car$ equals the coefficient of 
$t^n$ in $\pi(\ca,t)$.
\end{lemma}
\begin{proof}
It follows from the definition of general position that 
$L(\car)=L(\ca)_{\leq n-1}$, where the latter is
the truncation of the lattice $L(\ca)$ to rank $n-1$.  Therefore
$\pi(\ca,t)=\pi(\car,t)+bt^n$ for some $b$.  By a theorem of
Zaslavsky~\cite{zas75}, the number of bounded chambers of any arrangement
${\mathcal B}$ equals $(-1)^{{\rm rank}\,\mathcal B}\pi({\mathcal B},-1)$.
Substituting $t=-1$ shows $b$ is the number of
bounded chambers in $\car$, since $\ca$ itself has none.
\end{proof}
Let $\epsilon>0$ be a fixed choice of positive, real number.
\begin{lemma}\label{lem:gp}
For any $x\in C_\ba\cap(\R^n-V)$ 
let $H_x$ be the hyperplane normal to $x$, passing through $\epsilon x$.
Then $H_x$ is in general position to $\ca$.
\end{lemma}
\begin{proof}
Suppose $X\cap H_x=\emptyset$ for some nonzero intersection of hyperplanes
$X$.  Say $X=\cap_{\alpha\in S}H_\alpha$, where $S\subseteq\Phi$.  Since
$X\neq0$, the roots $S$ do not span $\R^n$.  Since $X$ and $H_x$ are parallel,
$x$ is a linear combination of the roots $S$; then $x\in V$, a contradiction.
\end{proof}
For each $y\in\R^n$ with $(x,y)>0$, let $y^{H_x}$ denote the unique, positive
multiple of $y$ which lies in $H_x$.  
Note that each chamber of $\cax$ has the form $C\cap H_x$ for some 
chamber $C$ of $\ca$.  If $C\cap H_x$ is bounded, then $C$ is just a cone
over $C\cap H_x$ with retraction $y\mapsto y^{H_x}$.  In particular,
$(x,y)>0$ for all $y\in C$.
For any $x\in\R^n-V$, let
\begin{equation}\label{def:bx}
B_x=\set{g\in W:\hbox{$(x,gx)>0$ and $(gx)^{H_x}$ is in a bounded 
chamber of $\cax$}}.
\end{equation}
Since $x\not\in V$, the orbit $Wx$ has exactly one point in each chamber
of $\ca$.  It follows that $\abs{B_x}$ is the number of bounded chambers
of $\cax$.
\begin{lemma}\label{lem:bij}
For any $x\in\R^n-V$, we have
$$
B_x=\set{g\in W:g^{-1}x\in\sigma_\ba}.
$$
\end{lemma}
\begin{proof}
A chamber $C\cap H_x$ of $\cax$
is bounded if and only if $C$ does not contain a ray in $H_x$.
Equivalently, all points in $C\cap H_x$ (or, just as well, in $C$)
have positive inner product with respect to $x$.  

That is, $g\in B_x$ if and only if, for all $y\in C_\ba$, 
$$
(gy,x)>0 \quad\Longleftrightarrow\quad (y,g^{-1}x)>0\quad
\Longleftrightarrow\quad g^{-1}x\in\sigma_\ba,
$$
since $\sigma_\ba$ is the normal cone to $C_\ba$.
\end{proof}

\subsection{Proof of Theorem~\ref{th:count}}
Fix a point $x\in\R^n-V$.  By construction, $x$ lies in some
(open) chamber $C$.  Without loss of generality, $C=C_\ba$.  
Let $H_x$ be the hyperplane normal to $x$, containing $\epsilon x$.
Using 
Lemmas~\ref{lem:trunc}, \ref{lem:gp}, and equation \eqref{eq:terao},
we see the number of bounded chambers in $\cax$ equals
$\prod_{i=1}^n(d_i-1)$.

On the other hand, the number of bounded chambers of $\cax$ equals
$\abs{B_x}$; by Lemma~\ref{lem:bij}, this equals the number of
$g\in W$ for which $x\in g\sigma_\ba$.\qed
\begin{figure}[htb]
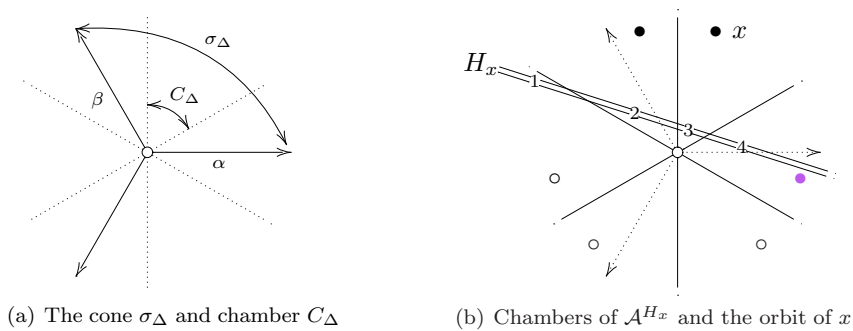

\begin{center}
\subfigure[The cone $\sigma_\ba$ and chamber $C_\ba$]{
\label{fig:onea}
\begin{minipage}[t]{0.45\linewidth}
\xygraph{
*{\cir<2pt>{}}="o"-@{>}_{\alpha}[]!{(1.5,0)}{}="A"
"o"-@{>}^{\beta}[]!{(-0.75,1.3)}{}="B"
"o"-@{>}[]!{(-0.75,-1.3)}{}
"o"-@{.}[]!{(0,1.5)}{}
"o"-@{.}[]!{(0,-1.5)}{}
"o"-@{.}[]!{(1.3,0.75)}{}
"o"-@{.}[]!{(-1.3,-0.75)}{}
"o"-@{.}[]!{(1.3,-0.75)}{}
"o"-@{.}[]!{(-1.3,0.75)}{}
[]!{(0,0.5)}-@{<->}@/^0.6ex/^{C_\ba}[]!{(0.433,0.25)}
"B"-@{<->}@`{(0.85,1.4)}^{\sigma_\ba}"A"[]!{!(0.35,0)}{}
}
\end{minipage}
}
\hskip4em
\subfigure[Chambers of $\cax$ and the orbit of $x$]{\label{fig:oneb}
\begin{minipage}[t]{0.45\linewidth}
\xy
%\UseCrayolaColors
%\newxyColor{unbd}{0.2}{grey}{}
\xygraph{
*{\cir<2pt>{}}="o"-@{.>}[]!{(1.5,0)}{}%{\alpha}="A"
"o"-@{.>}[]!{(-0.75,1.3)}{}[]!{!(-0.16,-0.04)}{}%{\beta}="B"
"o"-@{.>}[]!{(-0.75,-1.3)}{}
"o"-@{-}[]!{(0,1.5)}{}
"o"-@{-}[]!{(0,-1.5)}{}
"o"-@{-}[]!{(1.3,0.75)}{}
"o"-@{-}[]!{(-1.3,-0.75)}{}
"o"-@{-}[]!{(1.3,-0.75)}{}
"o"-@{-}[]!{(-1.56,0.9)}{}
[]!{(-1.5,0.75)}*{{\scriptstyle1}}="1"
-@{=}[]!{!(1.05,-0.333)}*{{\scriptstyle2}}
-@{=}[]!{!(.562,-0.179)}*{{\scriptstyle3}}
-@{=}[]!{!(.562,-0.179)}*{{\scriptstyle4}}
-@{=}[]!{(1.65,-0.25)}{}
"1"-@{=}[]!{!(-.562,0.179)}*{H_x}
[]!{(0.4,1.26)}*{\bullet}[]!{!(0.25,0)}*{x}
[]!{(-0.4,1.26)}*{\bullet}
[]!{(-1.29,-0.284)}*{\circ}
%{\scriptstyle+}
[]!{(1.29,-0.284)}*[unbd]{\bullet}
[]!{(-0.881,-0.976)}*{\circ}
[]!{(0.881,-0.976)}*{\circ}
}
\endxy
\end{minipage}
}
\end{center}
\caption{The $A_2$ root system}
\end{figure}
\begin{ex}
Let $\ba=\set{\alpha,\beta}$ be the base of the $A_2$ root system,
shown in Figure~\ref{fig:onea}.  Recall $d_1=2, d_2=3$; then 
$\nu(\sigma_\ba)=\frac{1\cdot2}{2\cdot3}$.  In Figure~\ref{fig:oneb},
the chambers of $\cax$ are labelled
$1$ through $4$.  As expected, two chambers (labelled 
$2$ and $3$) are bounded.  For a given $x\in C_\ba$, points
$gx$ in its orbit are
marked with a ``$\circ$'' if $(x,gx)\leq0$.  If $(x,gx)>0$, the 
point $gx$ is black where the chamber $(gx)^{H_x}$ is bounded and 
``$\xy *[unbd]{\bullet}\endxy$'' otherwise.
\end{ex}

\section{The identity}\label{sec:id}
Now suppose that $\Phi\subseteq\R^n$ is an irreducible,
crystallographic root system of
rank $n$.  Let $\widetilde{\Phi}$ denote the affine root system of $\Phi$,
with base $\widetilde{\ba}=\ba\cup\set{\alpha_0}$.  Let $\widetilde{D}$ 
denote the extended Dynkin diagram of $\Phi$.  For each simple root
$\alpha_i\in \widetilde{\ba}$, 
let $\Phi_i$ be the sub-root system of $\Phi$ with base
$\Delta_i=\widetilde{\ba}-\set{\alpha_i}$.  Then $\Phi=\Phi_0$, and recall that
the Dynkin diagram of $\Phi_i$ is obtained by deleting the vertex
corresponding to $\alpha_i$ from $\widetilde{D}$.

For each $i$, $0\leq i\leq n$, let $(d^{(i)}_1,\ldots,d^{(i)}_n)$ denote
the degrees of $\Phi_i$.  De Concini and Procesi found that, for each
irreducible type, an unexpected identity held:
\begin{theorem}[Theorem~1.2 of \cite{dCP06}]
For an irreducible, crystallographic root system $\Phi$ of rank $n$, 
\begin{equation}\label{eq:curious}
\sum_{i=0}^n\prod_{j=1}^n\frac{d_j^{(i)}-1}{d_j^{(i)}}=1.
\end{equation}
\end{theorem}
By (re)deriving their result from Theorem~\ref{th:vol}, a geometric
interpretation becomes apparent.
\begin{proof}
Let $A_0$ denote the fundamental alcove of $\Phi$.  This is a simplex
bounded by the (affine) reflecting 
hyperplanes $\set{H_{\alpha_i}\colon 0\leq i\leq n}$.  For each $i$,
let $v_i$ be the vertex of $A_0$ that is
opposite the face contained in $H_{\alpha_i}$.  The normal cone to $A_0$ at
$v_i$ is spanned by the vectors $\widetilde{\ba}-\set{\alpha_i}$, so
it is just the cone $\sigma_{\ba_i}$.
Then 
$$
\nu(\sigma_{\ba_i})=
\prod_{j=1}^n\frac{d_j^{(i)}-1}{d_j^{(i)}},
$$
by the volume formula
\eqref{eq:volformula}.
However, the normal cones to the vertices of 
any polytope partition a dense open subset of $\R^n$, so their volumes 
sum to $1$.  
\end{proof}
\begin{remark}
We have seen that the volume formula \eqref{eq:volformula} also holds for
finite, noncrystallographic root systems.  For the irreducible types,
\eqref{eq:volformula} gives
$$
\renewcommand{\arraystretch}{1.2}
\begin{array}{|l||l|l|l|}\hline
\hbox{Type} & I_2(m) & H_3 & H_4\\ \hline
\nu(\sigma_\ba) & (m-1)/(2m) & 3/8 & 6061/14\,400 \\ \hline
\end{array}
$$
Although the identity \eqref{eq:curious} no longer makes sense, one
might still be tempted to compute the left side formally for diagrams
that extend $H_3$ or $H_4$ by a vertex 
in such a way that all proper subdiagrams
are of finite type.  (These include the Coxeter groups $H_3^{\rm aff}$
and $H_4^{\rm aff}$ of Patera and Twarock,~\cite{PaTw02}.)
Perhaps unsurprisingly, however, 
an exhaustive search shows that the identity fails to hold 
for any such diagram.
\end{remark}

\providecommand{\bysame}{\leavevmode\hbox to3em{\hrulefill}\thinspace}
\providecommand{\MR}{\relax\ifhmode\unskip\space\fi MR }
\providecommand{\MRhref}[2]{%
  \href{http://www.ams.org/mathscinet-getitem?mr=#1}{#2}
}
\providecommand{\href}[2]{#2}

\end{document}